\documentclass[11pt]{amsart}

\usepackage{palatino}
\newcounter{EGnum}
\newcommand{\EGNumber}{\theEGnum\stepcounter{EGnum}}

\newenvironment{proofof}[1]{\noindent
               \textit{#1.}}{\hfill{{\sc $\square$}}\\}

\usepackage{latexsym} 
\usepackage{amsmath} 
\usepackage{amssymb} 
\usepackage{amscd} 
\usepackage{psfig} 
\usepackage{amsfonts}

\newtheorem{theorem}{Theorem}[section]

\newtheorem{lemma}[theorem]{Lemma} 
 
\newtheorem{proposition}[theorem]{Proposition} 
\newtheorem{definition}[theorem]{Definition} 
\newtheorem{remark}[theorem]{Remark}  

\newtheorem{DEF}[theorem]{Definition}

\input xy
\xyoption{all}

\def\to{\rightarrow} 
\def\into{\hookrightarrow}

\newcommand{\Z}{\mathbb Z}

\newcommand{\R}{\mathbb R}
\newcommand{\RR}{\mathcal R}
\newcommand{\C}{\mathbb C}
\newcommand{\CC}{\mathbb C}
\renewcommand{\L}{\mathbb L}
\newcommand{\LL}{\mathbb L}
\renewcommand{\i}{\ensuremath{\iota}}
\renewcommand{\P}{\mathbb P}
\newcommand{\PP}{\mathbb P}

\newcommand{\algk}{\ensuremath{\mathfrak{k}}}
\newcommand{\algt}{\ensuremath{\mathfrak{t}}}

\newcommand{\algh}{\ensuremath{\mathfrak{h}}}

\newcommand{\Ker}{\mbox{Ker}}
\newcommand{\Lie}{\mbox{Lie}}

\newcommand{\st}[1]{\ensuremath{^{\scriptstyle \textrm{#1}}}}

\numberwithin{equation}{section}

\def\iso{\cong}

\newcommand{\alphaparenlist}{
  \renewcommand{\theenumi}{\alph{enumi}}%
  \renewcommand{\labelenumi}{(\theenumi)}%
}

\newcommand{\Alist}{
  \renewcommand{\theenumi}{\arabic{enumi}}%
  \renewcommand{\labelenumi}{(A\theenumi)}%
}

\title[GKM theory for non-isolated fixed points]{GKM theory for torus actions with non-isolated fixed points}
\author{Victor Guillemin}
\address[Victor Guillemin]{Department of Mathematics\\Massachusetts Institute of Technology 2-251\\Cambridge MA 02139}
\email{vwg@math.mit.edu}
\author{Tara S.\ Holm}
\address[Tara Holm]{Department of Mathematics 3840\\University of California, Berkeley\\Berkeley, CA 94720-3840}
\email{tsh@math.berkeley.edu}

\begin{document}

\begin{abstract}
Let $M^{2d}$ be a compact symplectic manifold and $T$ a compact
$n$-dimensional torus.  A Hamiltonian action, $\tau$, of $T$ on
$M$ is a GKM action if, for every $p \in M^T$, the isotropy
representation of $T$ on $T_pM$ has pair-wise linearly
independent weights.  For such an action the projection of the set
of zero and one-dimensional orbits onto $M/T$ is a regular
$d$-valent graph; and Goresky, Kottwitz and MacPherson have
proved that the equivariant cohomology of $M$ can be computed
from the combinatorics of this graph.  (See \cite{GKM:eqcohom}.)  In
this paper we define a ``GKM action with non-isolated fixed
points'' to be an action, $\tau$, of $T$ on $M$ with the
property that for every connected component, $F$ of $M^T$ and $
p \in F$ the isotropy representation of $T$ on the normal space
to $F$ at $p$ has pair-wise linearly independent weights.  For
such an action, we show that all components of $M^T$ are
diffeomorphic and prove an analogue of the theorem above.
\end{abstract}

\maketitle

\section{Introduction}
\label{se:intro}

Let $M^{2d}$ be a compact symplectic manifold, $T$ an
$n$-dimensional torus and $\tau$ a Hamiltonian action of $T$ on
$M$.  We will denote by $M^T$ the fixed point set of $\tau$ and
by $M_1$ the set
\begin{displaymath}
  \{ x \in M \, , \, \dim T \cdot x =1 \} \, ;
\end{displaymath}
and we will say that $\tau$ is a GKM action if
\begin{equation}
  \label{eq:1.1}
  \# M^T < \infty
\end{equation}
and
\begin{equation}
  \label{eq:1.2}
  \dim M_1 =2 \, .
\end{equation}
Let
\begin{displaymath}
  V=\{ p_1 , \ldots , p_e \}
\end{displaymath}
be the elements of $M^T$ and
\begin{displaymath}
  E= \{ e_1 , \ldots , e_N \}
\end{displaymath}
the connected components of $M_1$. For each $e_i$ let
$\bar{e}_i$ be its closure in $M$.  From the assumptions
(\ref{eq:1.1}) and (\ref{eq:1.2}) one easily deduces:
\begin{enumerate}
\item Each $\bar{e}_i$ is an embedded copy of $\CC P^1$

\item The set  $\bar{e}_i -e_i$ is a two element subset of
 $V$.

\item For $i \neq j$,  $\bar{e}_i \cap \bar{e}_j$
is empty or is a one element subset of  $V$.
\item For $p \in V$, the set $\{ e_i , \, p
\in \bar{e}_i \}$ is a $d$-element subset of $E$.
\end{enumerate}
In other words $V$ and $E$ are the vertices and edges of a
regular $d$-valent graph, $\Gamma$.

Moreover, by noting how $T$ acts on each of the $\bar{e}_i$'s one
gets a \emph{labeling} of the oriented edges of $\Gamma$ by
elements of the weight lattice of $T$.  Explicitly, if $e$ is an
oriented edge of $\Gamma$ joining $p$ to $q$ one can assign to
$e$ the weight of the isotropy representation of $T$ on the
tangent space to $\bar{e}$ at $p$.  Denoting this weight by
$\alpha_e$, the assignment, $e \to \alpha_e$, defines a labeling
of the type above. We will call this labeling the
\emph{axial function} on the graph $\Gamma$.

We will show in \S~\ref{se:actions} how to construct from the data $(\Gamma
,\alpha)$ a commutative ring, $H(\Gamma ,\alpha)$, and sketch a
proof of the main result of \cite{GKM:eqcohom} which asserts that
\begin{equation}
  \label{eq:1.4}
  H (\Gamma ,\alpha) \simeq H_T (M;\C) = H_T^*(M) \, .
\end{equation}

Our main result is a generalization of the GKM theorem.  To describe it,
we recall that the hypotheses
(\ref{eq:1.1}) and (\ref{eq:1.2}) can be reformulated somewhat
differently.

\begin{proposition}
  \label{prop:1.1}
The conditions (\ref{eq:1.1}) and (\ref{eq:1.2}) are satisfied if
and only if, for every $p \in M^T$, the weights, $\alpha_{i,p}$,
$i=1,\ldots ,d$, of the isotropy representation of $T$ on $T_pM$
are pair-wise linearly independent, i.e.,~for $i \neq j$
$\alpha_{i,p}$ is not a multiple of $\alpha_{j,p}$.
\end{proposition}

This suggests imposing a slightly weaker ``GKM hypothesis'' on
the action $\tau$:

\begin{definition}
\label{def:1.2}
The action, $\tau$, is a {\em GKM action with non-isolated fixed
points} if, for every connected component, $F$ of $M^T$ and $p
\in F$ the isotropy representation of $T$ on the normal space to
$F$ at $p$ has pair-wise linearly independent weights.

\end{definition}

This relatively innocuous assumption has some surprising
implications.  Let
\begin{equation}
  \label{eq:1.5}
  \{ F_i \, ; \, \, i=1,\ldots ,\ell \}
\end{equation}
be the connected components of $M^T$, and let
\begin{equation}
  \label{eq:1.6}
  \{ W^0_i \, ; \, \, i=1,\ldots ,N \}
\end{equation}
be the connected components of $M_1$.

\begin{theorem}
  \label{th:1}

The sets above have the following properties:

\alphaparenlist
\begin{enumerate}
\item 
The $F_{i}$'s are all diffeomorphic and, in particular, are all
of the same dimension, $2m$.

\item 
  The closure, $W_i$, of $W^0_i$ is a symplectic submanifold of $
 M$ of dimension $2m+2$.

\item 
  $W^T_i$ is the union of two $X_j$'s.

\item 
$W_i \cap W_j$ is either empty or is a single $X_k$.

\item 
For $X=X_j$
\begin{displaymath}
  \# \{ i, \,\, X \subseteq W_i \} =d-m \, .
\end{displaymath}

\end{enumerate}

\end{theorem}

Thus, as above, the sets (\ref{eq:1.5}) and (\ref{eq:1.6}) are the vertices, $V$ and the
edges, $E$ of a regular $(d-m)$-valent graph, $\Gamma$.
Moreover, as above, this is a labeled graph.  If $W$ is one if
the $W_i$'s and $X$ one of the two connected components of $W^T$,
the isotropy representation of $T$ on the normal space to $X$ in
$W$ at $p \in X$ does not depend on $p$.  The weight of this
representation gives one a labeling, $e \to \alpha_e$, of the
oriented edges of $\Gamma$ by elements of the weight lattice of
$T$.  We will prove in \S\ref{se:gkm} the following generalization of (\ref{eq:1.4}).

\begin{theorem}
  \label{th:2}\label{th:MAIN}
For $F=F_i$,
\begin{equation}
  \label{eq:1.7}
  H_T (M) \simeq H (\Gamma ,\alpha) \otimes H(F) \, .
\end{equation}

\end{theorem}

The $W_i$'s, unlike the $X_i$'s, are not all diffeomorphic.
However, if $W$ is one of the $W_i$'s and $F$ one of the two
connected components of $W^T$, we will prove:

\begin{theorem}
  \label{th:3}
The normal bundle, $\LL$, of $X$ in $W$ is a complex line bundle
and
\begin{equation}
  \label{eq:1.8}
  W \simeq \PP (\LL \oplus \CC) \, .
\end{equation}

\end{theorem}

We will conclude this introduction with a brief summary of the
contents of this article.  In \S\ref{se:actions} we will discuss in more detail
the graph theoretic aspects of GKM theory.  In particular we will
describe what we mean by an ``action'' of a torus on a graph,
define the ``equivariant cohomology ring'' of a graph, and sketch
a proof of (\ref{eq:1.4}).

In \S\ref{se:fixedpts} we will prove theorems \ref{th:1} and \ref{th:3}, and in
\S\ref{se:gkm} we will prove theorem~\ref{th:2}.  In \S\ref{se:example} we will discuss a
few examples of ``GKM actions with non-isolated fixed points''.
All these examples are fiber bundles
\begin{displaymath}
  X \hookrightarrow M \overset{\pi}{\longrightarrow} F
\end{displaymath}
with the property that $T$ acts fiberwise and that the action
on the  fiber, $X$, is a GKM action in the usual sense.  
In \S\ref{se:sympl} we discuss the symplectic
structure of the $X_i$'s and show that the question of when two
$X_i$'s are \emph{sympelctomorphic} is closely related to the
question of when two $W_i$'s are \emph{diffeomorphic}.
Finally, in \S\ref{se:holonomy}, we discuss some holonomy invariants
of $M$ whose vanishing may imply that $M$ is a fiber product of the
type above.

\section{$T$ actions on graphs}
\label{se:actions}

Let $\Gamma$ be a regular $d$-valent graph and let $V$ and $E$ be
the vertices and the oriented edges of $\Gamma$.  For every $e
\in E$ we will denote by $\bar{e}$, the edge $e$, wth its
orientation reversed and we will denote by $i(e)$ and $t(e)$ the
intial and terminal vertices of $e$.  Thus $t(\bar{e})=i(e)$ and
$i(\bar{e}) =t(e)$.
Let $\RR_k (T)$ be the set of (equivalence classes of)
$k$-dimensional representations of $T$.  

\begin{definition}\label{def:action}
We define an
\emph{action of $T$ on} $\Gamma$ to be a pair of maps
\begin{displaymath}
  \tau : V \to \RR_d (T)
\end{displaymath}
and %
\begin{displaymath}
  \gamma :E \to \RR_1 (T)
\end{displaymath}
satisfying

\Alist
\begin{enumerate}
\item $\tau_p = \bigoplus_{i(e)=p} \gamma_e$ \, , 

\item $\gamma_{\bar{e}} = \gamma^*_e$ \, , and

\item $\tau_p (g) = \tau_q (g)$ for 
$p=i(e) \, , \, q=t(e)$ and $ g\in \Ker( \gamma_e)$.
\end{enumerate}
\end{definition}

\begin{remark}
By (A1), the mapping $\tau$ is determined by the mapping
$\gamma$.  Moreover, if for every oriented edge, $e$, we let
$\alpha_e$ be the weight of the representation $\gamma_e$,
$\gamma_e$ is determined by $\alpha_e$.  Hence $\gamma$ and
$\tau$ are determined by the axial function, $\alpha$,
which assigns to each oriented edge, $e$, the element,
$\alpha_e$, of the weight lattice of $T$.
\end{remark}

\begin{remark}
  The axioms (A1)--(A3) translate into axioms on $\alpha$.  For
  instance axiom~(A2) is equivalent to:
  $\alpha_{\bar{e}}=-\alpha_e$.
\end{remark}

\begin{remark}
Let $\tau$ be a GKM action of $T$ on $M$, and for each $p
  \in M^T$ let $\tau_p$ be the isotropy representation of $T$ on
  the tangent space to $M$ at $p$.  If $e$ is a connected
  component of $M_1$ we can orient $e$ by specifying that one of
  the two points in $\bar{e}-e$ is the initial vertex $i (e) =p$,
  of $e$ and the other the terminal vertex, $t(e) =q$.  If we let
  $\gamma_e$ be the isotropy representation of $T$ on the tangent
  space to $\bar{e}$ at $p$, the mappings, $p\to \tau_p$ ad $e
  \to \gamma_e$ define an action of $T$ on $\Gamma$.
\end{remark}

\begin{remark}
Let $\tau$ be a GKM action of $T$ on $M$ ``with
  non-isolated fixed points''.  Let $X_i$, $i=1, \ldots , \ell$,
  and $W^0_i$, $i=1,\ldots,N$, be the connected components of
  $M^T$ and $M_1$; and for each $W^0_i$, let $W_i$ be its closure
  in $M$.  Then, by Theorem~\ref{th:1} the sets
$$
    V = \{ X_i , \quad i=1,\ldots , \ell \}
$$
and
$$
    E = \{ W_i \, , \quad i=1,\ldots ,N \}
$$
are the vertices and edges of a regular $(d-m)$-valent graph.  For
each $X_i$, let $\tau_i$ be the isotropy representation at $p \in
X_i$ of $T$ on the normal space to $X_i$ at $p$.  (Since $X_i$ is
connected this representation doesn't depend on $p$.)  Similarly,
for every $W_i$ et $X_j$ and $X_k$ be the connected components of
$W_i-W^0_i$ and let $\gamma_{j,k}$ be the isotropy
representaiton, at $p \in X_j$, of $T$ on the normal space at $p$
to $X_j$ in $W_i$.  Then the mappings, $i \to \tau_i$ and
$(j,k)\to \gamma_{j,k}$ define an action of $T$ on $\Gamma$.
\end{remark}

Let $\algt = \Lie(T)$ and let $S(\algt^*)$ be the ring of polynomial
functions on $\algt$.  Given a graph $T$ and an action of $T$ on
$\Gamma$ we will deine the \emph{equivariant cohomology ring},
$H(\Gamma ,\alpha)$, of $\Gamma$ to be the set of maps
\begin{displaymath}
  f:V \to S(\algt^*)
\end{displaymath}
which, for all $e \in E$, satisfy the compatability condition
\begin{equation}
  \label{eq:2.1}
  f_p -f_q \in \alpha_e \cdot S (t^*)
\end{equation}
for $p=i(e)$ and $q=t(e)$.  The GKM theorem asserts

\begin{theorem}
  \label{th:2.1}

If $M$ is a GKM manifold
\begin{displaymath}
  H(\Gamma ,\alpha) \simeq H_T (M) \, .
\end{displaymath}
\end{theorem}

We will give a brief sketch of how to prove this since we will
prove Theorem~\ref{th:2} by mimicking this proof.

\subsection*{Step 1.\quad (The Kirwan formality theorem.)}

This asserts that as an $S(\algt^*)$ module,
\begin{equation}
  \label{eq:2.2}
  H_T (M) = H(M) \otimes S(\algt^*) \, .
\end{equation}
By a theorem of Borel, the restriction map
\begin{equation}
  \label{eq:2.3}
  r: H_T (M) \to H_T (M^T)
\end{equation}
has, as kernel, the torsion elements in $H_T (M)$; however,
Kirwan's theorem implies that $H_T(M)$ is a free $S(\algt^*)$ module
hence (\ref{eq:2.3}) is injective.  Moreover, since $T$ acts
trivially on $M^T$ and $M^T$ is finite, $H_T(M^T)$ is a sum of
copies of $S(\algt^*)$
\begin{equation}
  \label{eq:2.4}
  \bigoplus_{p\in M^T} S(\algt^*)_p 
\end{equation}
or alternatively is the set of maps,
\begin{displaymath}
  f:V \to S(\algt^*) \, .
\end{displaymath}
We claim

\begin{lemma}
  The image of $r$ is contained in the subring, $H(\Gamma
  ,\alpha)$ of the ring of maps of $V$ into $S(\algt^*)$.
\end{lemma}

\begin{proof}
Let $e^0$ be a connected component of $M_1$, let $e$ be its
closure and let $p$ and $q$ be the elements of $e-e^0$.  The
kernel $T_e$ of $\gamma_e : T \to S^1$ acts trivially on $e$;
therefore, denoting by $\algt_e$ the Lie algebra of $T_e$
\begin{displaymath}
  H_{T_e} (e) = H(e) \otimes S(\algt^*_e) \, .
\end{displaymath}
In particular, letting $i_p$ and $i_q$ be the inclusions of $p$ and
$q$ into $e$ the induced maps
\begin{eqnarray*}
  i^*_p : && H_{T_e} (e) \to S (\algt^*_e)\\
\noalign{\hbox{and}}\\
i^*_q : && H_{T_e} (e) \to S(\algt^*_e)
\end{eqnarray*}
are identical.  In particular for every $f$ in the image of the
restriciton map (\ref{eq:2.3}) $f_p$ and $f_q$ have to satisfy
the compatibility condition (\ref{eq:2.1}).

\subsection*{Step 2. \quad Betti numbers.}

Fix a vector $\xi \in t$ such that for all $e \in E$ $\alpha_e
(\xi) \neq 0$ and for every $p \in V$ let 
\begin{equation}
  \label{eq:2.5}
  \sigma_p = \# \{ \alpha_e \, ; 
     2(e) =p \hbox{  and  } \alpha_e (\xi)<0 \} \, .
\end{equation}
We define the $i$\st{th} Betti number $\beta_i (\Gamma)$ of the
graph $\Gamma$ to be:
\begin{equation}
  \label{eq:2.6}
  \# \{ p \in V \, , \quad \sigma_p =i \} \, .
\end{equation}
The numbers (\ref{eq:2.5}) depend on the choice of $\xi$, however
one can show that the numbers (\ref{eq:2.6}) don't.  Moreover,
one can prove by elementary Morse theory that $\beta_{2i+1}(M)
=0$ and
\begin{equation}
  \label{eq:2.7}
  \beta_i (\Gamma) = \beta_{2i}(M) \, .
\end{equation}
Thus by the Kirwan formality theorem
\begin{equation}
  \label{eq:2.8}
  \dim H^{2k}_T (M) = \sum \beta_i (\Gamma)
  \dim S (\algt^*)^{k-i} \, .
\end{equation}
Note that since the odd Betti numbers of $M$ are
zero, \eqref{eq:2.2} implies that the odd equivariant cohomology
groups of $M$ are zero.

\subsection*{Step 3. \quad Graph theoretic Morse inequalitites.}

These assert that
\begin{equation}
  \label{eq:2.9}
  \dim H^k (\Gamma ,\alpha) \leq \sum \beta_i (\Gamma )
  \dim S(\algt^*)^{k-\ell} \, .
\end{equation}
For the relatively elementary proof of these inequalities, see
\cite{GZ:graphs}.  Combining steps 1, 2 and 3 we conclude
that the map (\ref{eq:2.3}) maps $H_T(M)$ bijectively onto
$H(\Gamma ,\alpha)$.
\end{proof}

\section{The fixed points}\label{se:fixedpts}

Let $M$ be a compact connected symplectic manifold of dimension $2d$,
and let $T$ be an $n$-torus and let $\tau$ be a Hamiltonian action of
$T$ on $M$ with moment map $\Phi:M\to\algt^*$.  We make the following
GKM assumption: For every connected component $F$ of $M^T$, the
weights
\begin{equation}\label{eq:weights}
\alpha_{i,F} \mbox{\phantom{boo}} i=1,\dots,r
\end{equation}
of the isotropy represention of $T$ on the normal bundle to $F$ are
pairwise linearly independent.  As above, we call such a manifold a
non-isolated GKM manifold. We begin by examining the fixed point 
components of $M^T$, in the case when $T=S^1$ is a circle.

\begin{lemma}\label{le:codim}
Let $S^1$ act on a compact, connected symplectic manifold $M$ in
Hamiltonian fashion, with moment map $\phi:M\to\R$.  Every connected
component of $M^{S^1}$ is of codimension $2$ in $M$. 
\end{lemma}

\begin{proof}
Let $E$ be such a component and let $N\to E$ be the normal bundle to
$E$ in $M$. The
weights of the isotropy representation of $T$ on $N^H$ all have to be
multiples of $\alpha$, since $\algt^*$ is one-dimensional. So, because
the weights must be two-independent by the GKM assumption above, there
can be only one weight, and so $\dim_{\R}(N^H)=2$.
\end{proof}

At any point $p\in F$, the Darboux theorem for $\phi$ says that there exists a
Darboux co\"ordinate system centered at $p$: co\"ordinates
$x_1,y_1,\dots,x_d,y_d$ such that locally near $p$, $F$ is defined by
$x_1=y_1=0$ and 
$$
\phi=\phi(p)+\alpha(\xi)(x_1^2+y_1^2).
$$
Thus, the component $F$ is either a maximum or a minimum of $\phi$ depending on
whether $\alpha(\xi)<0$ or $\alpha(\xi)>0$.  The same is true for
every other component of $M^{S^1}$.  However, by the Atiyah convexity
theorem, $\phi$ has at most one connected level set where it takes its
maximum value and one connected level set where it takes its minimum
value.  Thus, $M^{S^1}$ has exactly two connected components, $F$ and
$E$.

\begin{lemma}\label{le:diff}
The components $F$ and $E$ are diffeomorphic.
\end{lemma}

\begin{proof}
We can assume $\phi=0$ on $F$ and $\phi=1$ on $E$.  Let $\L$ be the
normal bundle to $F$ in $X$.  We can regard $\L$ as a complex line
bundle.  By the equivariant tubular neighborhood theorem, the action
of $T$ on $X$ is identical, near $F$ with the linear action of $T$ on
$\L$ and $\phi$ is just the length-squared function for a Hermitian
metric on $\L$.  Thus, for $c$ close to $0$,
$$
\phi^{-1}(c)/{S^1}=F.
$$
But all level sets, $\phi^{-1}(c)$ for $0<c<1$ are equivariantly
diffeomorphic, since there are no critical values between $0$ and $1$.
So all the reduced spaces
$$
\phi^{-1}(c)/S^1 \phantom{boo} 0<c<1
$$
are diffeomorphic.  But for $c$ close to $1$,
$$
\phi^{-1}(c)/{S^1}=E.
$$
by the same argument as above.  Hence, $F$ and $E$ are diffeomorphic.
\end{proof}

The result can be sharpened.  From the action of $T$ on $\L$ one gets
an action of $T$ on the bundle
$$
P(\L\oplus \C).
$$

\begin{theorem}\label{th:projBdle}
The component $M$ and $P(\L\oplus\C)$ are isomorphic as
$S^1$-manifolds. 
\end{theorem}

\begin{proof}
Equip $L$ with an inner product, and let $\psi : L \to [0,1)$ be the
function defined by 
$$
\psi (x,v) = \frac{|v|^2}{|v|^2 + 1}.
$$
This extends to a Morse-Bott function
$$
\psi : P(\L \oplus \C) \to [0,1]
$$
whose critical sets, $\psi^{-1}(1)$ and $\psi^{-1}(0)$, are copies of
$F$. The set $\phi^{-1}([\varepsilon , 1-\varepsilon])$ can be
identified with $F\times [\varepsilon, 1-\varepsilon])$.  So
$$
M_\varepsilon = \phi^{-1}([\varepsilon , 1-\varepsilon]) \to F\times
[\varepsilon, 1-\varepsilon])
$$
is a circle bundle with Chern class $c(\L)$.  Similarly,
$$
P^\varepsilon(\L\oplus\C) = \psi^{-1}([\varepsilon , 1-\varepsilon])
$$
is a circle bundle over this set with Chern class $c(\L)$.  Hence,
these bundles are isomorphic as manifolds with boundary.  One obtains
$P(\L\oplus \C)$ and $M$ from these manifolds by collapsing the circle
orbits on the boundary to points; however, the boundaries are
isomorphic as circle bundles.  Thus, the spaces obtained by this
collapsing are isomorphic.
\end{proof}

Now let $T$ be an $n$-torus, for $n>1$, and $\tau$ a GKM action on $M$
with non-isolated fixed points.

\begin{theorem}
All connected components of $M^T$ are diffeomorphic. 
\end{theorem}

\begin{proof}
Let $F$ be a connected component of $M^T$.  If $\alpha$ is one of the
isotropy weights \ref{eq:weights} on the normal bundle to $F$, let
$$
\algh = \{\xi\in\algt\ |\ \alpha(\xi)=0\},
$$
and let $H$ be the subtorus of $T$ with $Lie(H)=\algh$.  Let $X$ be
the connected component of $M^H$ containing $F$.

Now let $\xi\in\algt\setminus\algh$ and let $\phi^\xi$ be the restriction
of $\Phi^{\xi}$ to $X$.  We can regard $\phi^\xi$ as the moment map
associated with the action of the circle $T/H$ on $X$.  By
Lemma~\ref{le:diff},  $X$
contains precisely two components of $M^T$, and they are
diffeomorphic.  Let $\Gamma$ be the GKM
graph associated with the $T$-action on $M$.  To show that all
connected components of $M$ are diffeomorphic, we must show that
$\Gamma$ is connected.  Let $\phi$ be a Bott-Morse component of $\Phi$
with
$$
Crit(\phi)=M^T.
$$
That is, $\phi$ is a generic component of $\Phi$.  From $\phi$,
$\Gamma$ inherits a poset structure, and for each connected component
$\Gamma_0$ of $\Gamma$, the vertex of $\Gamma_0$ at which $\phi$ takes
its minimum corresponds to a component of $M^T$ at which $\phi$ takes
on a minimum value.  But the Atiyah convexity theorem says that there
is a unique component of $M^T$ on which this can occur.  Thus,
$\Gamma$ is connected, completing the proof of the theorem.
\end{proof}

\begin{remark}
The generic component $\phi$ of $\Phi$ in the proof above is a perfect
Bott-Morse function, with diffeomorphic critical sets.  We will use
this fact below to analyze the equivariant topology of $M$.
\end{remark}

\begin{remark}
The connected components of $M^T$ are not symplectomorphic.  We will
further discuss these symplectic structures in Section~\ref{se:sympl}.
\end{remark}

\section{GKM theory}\label{se:gkm}

Suppose $M$ is a compact, connected symplectic manifold, and
$\tau:M\times T \to M$ a
GKM action with non-isolated fixed points.  Let $\Gamma$ be the GKM graph
associated to $M$ and $\tau$.  Each edge $e$ of $\Gamma=(V,E)$ is labeled by a
one-dimensional representation $\gamma_e$ of $T$. Let $T_e$ be the
kernel of $\gamma_e$, and $\algt_e$ its Lie algebra.  To each vertex
of the graph,
$\Gamma$, we attach the ring $R=H_T^*(F)=H^*(F)\otimes S(\algt^*)$,
and to each edge, $e$, the
ring $R_e=H_{T_e}^*(F)=H^*(F)\otimes S(\algt_e^*)$. The map
$$
S(\algt^*)\to S(\algt_e^*)
$$
induces a map
$$
\pi_e:R\to R_e.
$$

\begin{DEF}
A map $f:V\to R$ is {\em $\Gamma$-compatible} if, for every edge $e=(p,q)$,
\begin{equation}\label{eq:*}
\pi_e(f(p))=\pi_e(f(q)).
\end{equation}
Let $H^*(\Gamma,F)$ be the ring of all of these maps.
\end{DEF}

We recall that Theorem~\ref{th:MAIN} asserts 
$$
H_T^*(M)=H^*(\Gamma,F).
$$
To prove this, we note that the Kirwan map
$$
H_T^*(M)\into H_T^*(M^T)
$$
is an injection.  Hence, if $F_\ell$ are the connected
components of $M^T$,  for $\ell=1,\dots,N$, then $H_T^*(M)$ sits inside the ring
$$
H_T^*(M^T)=\bigoplus_\ell H_T^*(F_\ell)
$$
and since $F_i\iso F$, $H_T^*(F_i)=H_T^*(F)=R$.  Thus, $H_T^*(M^T)$
consists of $N$ copies of $R$, each labeled by a vertex of $\Gamma$.
In other words,
$$
H_T^*(M^T)=Maps(V,R),
$$
and in particular, $H^*(\Gamma,F)$ is a subring of $H_T^*(M^T)$.  We
now prove the main theorem.  We break the proof of the theorem into
three steps.

\medskip

\begin{proofof}{Proof of Theorem~\ref{th:MAIN}}

\noindent {\bf Step One}: We will show that the restriction map
\begin{equation}\label{eq:1}
H_T^*(M)\into H_T^*(M^T)
\end{equation}
maps $H_T^*(M)$ into $H^*(\Gamma,F)$.

Let $\L$ be a line bundle over $F$, and let $X=P(\L\oplus \C)$.  Let
$\tilde{\L}$ be the tautology bundle over $X$, and let $x$ be its
Chern class.  By the Leray-Hirsch theorem,
\begin{equation}\label{eq:2}
H^*(X)=H^*(F)\oplus x\cdot H^*(F).
\end{equation}
Let $\i_+:F\to X$ be the embedding of $F$ onto $P(\{ 0\}\oplus\C)$.
Then $\i_+^*\tilde{\L}=1$, so $i_+^*x=0$.  Similarly, identifying
$P(\L\oplus\C)$ with $P(\C\oplus\L^{-1})$, and letting $\i_-:F\to X$
be the embedding onto $P(\C\oplus\{ 0\})$, one has $\i_-^*x=0$.  Thus,
by \eqref{eq:2},
\begin{equation}\label{eq:3}
\i_+^*c=\i_-^*c
\end{equation}
for every cohomology class $c\in H^*(X)$.

Now let $F=F_i$, let $T_e$ be a codimension one subgroup of $T$ and
let $X$ be a
connected component of $M^{T_e}$ containing $F$.  Let $\i_\pm:F\to X$
be the embeddings onto the two components of $X^T$. Then by
\eqref{eq:3}, the maps
$$
\xymatrix{
H_T^*(X)\ar[r] &  H_{T_e}^*(X)=H^*(X)\otimes
S(\algt_e^*)\ar[r]^{\hspace{0.7in}\i_\pm^*} & R_e
}
$$
are identical.  Thus, if $f$ is in the image of \eqref{eq:1}, it has
to satisfy compatibility conditions \eqref{eq:*} for every edge $e$ of
the graph $\Gamma$.

\medskip
\noindent {\bf Step Two}:
We first notice that $H^*(\Gamma,F)=H^*(\Gamma)\otimes H^*(F)$ and
$R=H^*(F)\otimes S(\algt^*)$. Hence, tensoring \eqref{eq:2.9} with $H^*(F)$,
we get the inequality
\begin{equation}\label{eq:6}
\dim H^k(\Gamma,F)\leq \sum_i\dim R^{k-i}\beta_i,
\end{equation}
where the $\beta_i$ are the Betti numbers of $\Gamma$.

\medskip
\noindent {\bf Step Three}:  Recall that $F_\ell$ is a connected
component of the critical set of the Bott-Morse function
$\phi=\Phi^\xi$.  Using Bott-Morse theory, Atiyah proves
\begin{equation}\label{eq:7}
\dim H_T^k(M)=\sum \dim H_T^{k-d_\ell}(F_\ell),
\end{equation}
where $d_\ell$ is the index of $F_\ell$ (see \cite{atiyah}).  Moreover, by elementary
Morse theory, it is easy to show that
\begin{equation}\label{eq:8}
\frac{d_\ell}{2}= \sigma_p
\end{equation}
where $p$ is the vertex of $\Gamma$ corresponding to $F_\ell$.
Thus, from \eqref{eq:7} and \eqref{eq:8}, one deduces
\begin{equation}\label{eq:9}
\dim H_T^*(M)=\sum_i\dim R^{k-i}\beta_i.
\end{equation}
Finally, the identiy \eqref{eq:9} and the inequality \eqref{eq:6}
imply that the Kirwan injection is a bijection of $H^*_T(M)$ onto
$H^*(\Gamma, F)$.
\end{proofof}

\section{Examples}\label{se:example}

We now describe several examples of GKM actions with non-isolated
fixed points, all of which start with the same 
basic set-up.  Let $F$ be a symplectic manifold, and let $\L_i$, for
$i=1,\dots,n$ be complex line bundles over $F$.  Let
$$
E=\L_1\oplus\cdots\oplus\L_n.
$$
The group $T=T^n$ acts on $E$  by acting fiberwise on the fiber
$$
E_p=(\L_1)_p\oplus\cdots\oplus(\L_n)_p
$$
by the action
$$
\tau(e^{i\theta})(v_1,\dots,v_n)=(e^{i\theta_1}\cdot v_1,\dots,
e^{i\theta_n}\cdot v_n).
$$
We use minimal coupling to produce a symplectic form $\omega$ on $E$,
and $T$ acts in a Hamiltonian fashion with respect to $\omega$.

\subsection{Projective bundles}
Let
\begin{equation}\label{eq:eg1}
\P(E)\to F
\end{equation}
be the projectivization of $E$.  From the action of $T$ on $E$, one
gets a fiber-wise action of $T$ on $\P(E)$ whose fixed points are
copies of $F$, and which satisfies our non-isolated GKM axiom.  The
graph associated to this space in the $n$-simplex, the same graph
associated to complex projective space in the ordinary GKM setting.

\subsection{Grassmannian bundles}
In a similar vein, let
\begin{equation}\label{eq:eg2}
\mathcal{G}r_k(E)\to F
\end{equation}
be the fiber bundle over $F$ whose fiber at $p$ is the Grassmannian
$\mathcal{G}r_k(E_p)$.  The $T$ action on $E$ defines a fiberwise
action of $T$ on the fiber of \eqref{eq:eg2} which also satisfies our
non-isolated GKM axiom. The graph associated to this space in the
Johnson graph $J(n,k)$, whose vertices consist of $k$-element subsets
of an $n$-element set. This is the same graph associated to the complex
Grassmannian $\mathcal{G}r(n,k)$ in the ordinary GKM setting.

\subsection{Partial flag bundles} One can continue in the
vein of \eqref{eq:eg1} and \eqref{eq:eg2} and take fiber bundles with
fiber some partial flag variety of $E_p$.

\subsection{Toric bundles}
In the last example, we apply symplectic reduction to the above
examples.  The torus $T^n$ acts fiberwise on the bundle $E$.  We may
make a symplectic reduction by $T^k$ to obtain
$$
\xymatrix{
\C^n/\!/ T^k\ar@{^{(}->}[r] & E/\!/ T^k \ar[d] \\
                           & F
}
$$
a fiber bundle over $F$ whose fiber is a complex toric variety.  In
this case, the moment image of $E$ is a simple convex polytope.

\medskip

These examples all exemplify the situation where there is a fiber
bundle $M\to F$ with fiber $X$, and a fiberwise action of $T$.  Modulo
assumptions on $F$ and $X$, the Leray-Hirsch theorem asserts
\begin{equation}\label{eq:eg3}
H_T^*(M)\iso H^*(F)\otimes H^*_T(X).
\end{equation}
Hence, if $X$ is a GKM manifold, one gets from \eqref{eq:eg3} the same
result as Theorem~\ref{th:MAIN}.

\section{The symplectic forms on the fixed point sets}\label{se:sympl}

Let $F$ be a connected component of $M^T$, $H$ a codimension one
subgroup of $T$, and $X$ a connected component of $M^H$ containing
$F$.  Let $F'$ be the other connected component of $M^T$ in $X$.
Without loss of generality, we may assume that $T/H$ acts faithfully
on $X$.  Let $\xi\in \algt\setminus\algh$ be the generator of this
group, normalized so that $\exp (2\pi\xi) = 1$, and let $\phi^\xi$ be
the $\xi$-component of the moment map $\Phi$.  Replacing $\xi$ by
$-\xi$ if necessary, we may assume that the restriction of $\phi^\xi$
takes its minimum value on $F$.

Now let $e$ denote the oriented edge of $\Gamma$ joining the vertices
corresponding to $F$ and $F'$.    We will call the difference
\begin{equation}
a_e = \phi^\xi(F') - \phi^\xi(F)
\end{equation}
the {\em length} of $e$.  Let $\L_e$ be the normal bundle to $F$ in
$X$.  If $\overline{e}$ is the edge $e$ with its orientation reversed,
then $\L_{\overline{e}}$ is the normal bundle to $F'$ in $X$.  So in
view of the isomorphisms
$$
X\iso \C P(\L\oplus 1) \iso \C P(1\oplus \L^{-1}),
$$
we have
\begin{equation}
\L_e\iso \L_{\overline{e}},
\end{equation}
and hence if $c_e$ and $c_{\overline{e}}$ are the Chern classes of
$\L_e$ and $\L_{\overline{e}}$ respectively, then
\begin{equation}\label{eq:6.3}
c_{\overline{e}} = - c_e.
\end{equation}
From the Duistermaat-Heckman theorem, one deduces the following theorem.

\begin{theorem}
Let $\omega_F$ and $\omega_{F'}$ be the symplectic forms on $F$ and
$F'$.  Then
\begin{equation}
[\omega_{F'}] = [\omega_F] + a_ec_e.
\end{equation}
\end{theorem}

We will discuss a few other properties of the assignment $e\mapsto
c_e$. We recall that a GKM action with non-isolated fixed points has
the following property: for every connected component $F$ of $M^T$,
and for every $p\in F$, the weights of the isotropy representation of
$T$ on the normal space to $F$ at $p$ are pairwise linearly
independent.  We will call a representation with this property
$2$-{\em independent}. More generally, we will call a representation of
$T$ with weights $\alpha_i$, for $i=1,\dots,N$ {\em $k$-independent} if
every subset of $k$ weights is linearly independent.  Now let $\Gamma$
be a regular $d$-valent graph, as in Section~\ref{se:actions}, and
$(\tau,\gamma)$ an action of $T$ on $\Gamma$.  For $p\in V$, let
$$
E_p = \{ e\in E\ |\ i(e) = p\}.
$$
A {\em connection} on $\Gamma$ is a function, $\nabla$, which assigns
to ever edge $e\in E$ with $p = i(e)$ and $q = t(e)$ a bijective map
$$
\nabla_e : E_p \to E_q
$$
with the property
$$
\nabla_{\overline{e}} = \nabla_e^{-1}.
$$
We say the connection is {\em compatible} with the action of $T$ if,
for every edge $e\in E$,
\begin{equation}
\nabla_e e' = e'' \Longrightarrow a_{e'} \cong a_{e''} \mod a_e.
\end{equation}
We leave the following strengthening of Axiom (A3) of
Definition~\ref{def:action} as an easy exerciese.

\begin{lemma}
Suppose that for all $p\in V$, $\tau_p$ is $3$-independent.  Then there
exists a unique $T$-compatible connection $\nabla$ on $\Gamma$.
\end{lemma}

In particular, let $\Gamma$ be the graph associated to $M$.  Then by
definition, $M$ is GKM if for every $p\in V$, $\tau_p$ is $2$-independent.
We make the stronger assumption that for every $p\in V$, $\tau_p$ is
$3$-independent.  Then the connection $\nabla$ is not only compatible with
the action of $T$, but is also compatible with the assignment
$e\mapsto c_e$.  More explicitly, let $e_1,\dots,e_k$ be the oriented
edges of $\Gamma$ with initial vertex $F$, and let $e_1',\dots,e_k'$,
with $e_1' = \overline{e_1}$, be the oriented edges of $\Gamma$ with
initial vertex $F'$.  We will order these edges so that $\nabla_{e_1}$
maps $e_i$ to $e_i'$.  We claim that
\begin{equation}\label{eq:6.6}
c_{e_i'} = c_{e_i} \mbox{ for } i = 2,\dots, k,
\end{equation}
and
\begin{equation}\label{eq:6.7}
c_{e_1} = -c_{e_1'}.
\end{equation}
\begin{proof}
The second identity \eqref{eq:6.7} follows immediately from
\eqref{eq:6.3}.  To prove \eqref{eq:6.6}, let $a_i = a_{e_i'}$, and
let $K_i = \exp(\algk_i)$, where 
$$
\algk_i = \{ \xi\in\algt\ |\ a_i(\xi) = a_1(\xi) = 0\}.
$$
Let $Y_i$ be the connected component of $M^{K_i}$ containing $X$.
Since the action of $T$ is $3$-independent, the codimension of $X$ in $Y_i$
is $2$.  Let $\L_i$ be the normal bundle of $X$ in $Y_i$.  Then the
restriction of $\L_i$ to $F$ is $\L_{e_i}$, and the restriction of
$\L_i$ to $F'$ is $\L_{e_i'}$.  Hence if $c_i$ is the Chern class of
$\L_i$ in $H^2(X;\Z)$, its restriction to $F$ is $c_{e_i}$, and its
restriction to $F'$ is $c_{e_i'}$.  Thus, by \eqref{eq:3},
$c_{e_i}=c_{e_i'}$.
\end{proof}

\section{Holonomy and diffeomorphisms of the fixed point sets}
\label{se:holonomy}

The GKM manifolds $M$ discussed in Section~\ref{se:example} are all twisted
products of $F$ with GKM manifolds with isolated fixed points.  We
now describe a holonomy invariant of the graph of $M$ which
measures the failure of $M$ to be such a twisted product.  Let $E$ and
$F$ be components of $M^T$.  Then, as we explain in
Section~\ref{se:fixedpts}, one can construct a diffeomorphism of $E$
onto $F$ using Morse 
theory.  This diffeomorphism is not unique, but it is easy to see that
it is unique up to isotopy.

In particular, suppose that $\gamma$ is a closed path in $\Gamma$
whose initial and terminal vertices are the vertex corresponding to
$F$.  Then by composing the diffeomorphisms above, one can associate
to $\gamma$ a diffeomorphism of $F$ onto itself which is unique up to
isotopy. Thus, letting $G=\pi_0(Dif(M))$ and letting
$\pi_1(\Gamma,F)$ be the fundamental group of $\Gamma$ with base point
$F\in V$, one gets a homomorphism
\begin{equation}\label{eq:7.1}
\Theta: \pi_1(\Gamma,F)\to G,
\end{equation}
which is the holonomy invariant alluded to above.

If $M$ is K\"{a}hler, and the action of $T$ preserves the K\"{a}hler
structure, one has a slightly more refined version of this invariant.
Namely, in this case, the action of $T^n$ on $M$ extends to a
holomorphic action of the complex torus $T^n_{\C}=P\C^*)^n$. In
particular, if $n=1$, one has a $\C^*$ action on $M$ and the
diffeomorphism between $E$ and $F$ in Lemma~\ref{le:diff} is given
explicitly by the map
$$
e\in F\mapsto y\in F
$$
if and only if, for some $m\in M$,
$$
x = \lim_{z\rightarrow 0} \tau_z(m)
$$
and
$$
y = \lim_{z\rightarrow \infty} \tau_z(m).
$$
In this case, the diffeomorphism is canonically defined and is a
biholomorphism.  Hence, in the K\"{a}hler case, the holonomy invariant
becomes a homomorphism
\begin{equation}\label{eq:7.2}
\Theta: \pi_1(\Gamma,F)\to Bihol(F).
\end{equation}
It is clear that if $M$ is a fiber product of the type discussed in
\S~\ref{se:example}, then these invariants vanish.  We do not know if
the converse is true.

\end{document}